\documentclass{amsart}
\usepackage{amsfonts}
\usepackage{amssymb}
\usepackage[latin1]{inputenc}
\newtheorem{thm}{Theorem}[section]

\newtheorem{lem}[thm]{Lemma}
\newtheorem{prop}[thm]{Proposition}
\theoremstyle{definition}

\theoremstyle{remark}

\numberwithin{equation}{section}

\newcommand{\boldf}{{\bf f}}

\newcommand{\cZ}{\mathcal Z}

\newcommand{\bbR}{\mathbb R}
\newcommand{\bbT}{\mathbb T}
\newcommand{\bbC}{\mathbb C}

\newcommand{\bbK}{\mathbb K}
\newcommand{\bbZ}{\mathbb Z}

\begin{document}

\title{Convergence versus integrability in Poincaré-Dulac normal form}

\author{Nguyen Tien Zung}
\address{GTA, UMR 5030 CNRS, Département de Mathématiques, Université Montpellier II}
\email{tienzung@math.univ-montp2.fr {\it URL}: www.math.univ-montp2.fr/\~{}tienzung}

\date{V.1, May 2001, V.2 Mars 2002, math.DS/0105193, submitted to Math. Res. Lett.}
\subjclass{37G05,70K45,34C14,70GXX}

\keywords{Poincaré-Dulac normal form, torus action, integrability, commuting vector
fields, isochore vector fields}%

\begin{abstract}
We show that, to find a Poincaré-Dulac normalization for a vector field is the same
as to find and linearize a torus action which preserves the vector field. Using this
toric characterization and other geometrical arguments, we prove that any local
analytic vector field which is integrable in the non-Hamiltonian sense admits a
local convergent Poincaré-Dulac normalization. These results generalize the main
results of our previous paper \cite{ZungBirkhoff} from the Hamiltonian case to the
non-Hamiltonian case. Similar results are presented for the case of isochore vector
fields.
\end{abstract}
\maketitle


\section{Introduction}

One of the classical problems in ordinary differential equations and dynamical
systems is to find necessary and sufficient conditions for the existence of a local
analytic (i.e. convergent) Poincaré-Dulac normalization for an analytic vector field
near an equilibrium point. Roughly speaking, there are (at least) two approaches to
this problem : analytical and geometrical. In the analytical approach, with
classical results due to Poincaré, Siegel, Bruno and others, one uses Diophantine
conditions to control small divisors which appear in a step-by-step normalization
process, and uses the fast converging iteration method to arrive at convergence
results, see e.g. \cite{Bruno,Roussarie}. In the geometrical approach, one replaces
Diophantine conditions by symmetries and first integrals, and replaces the fast
convergence method by arguments of geometrical nature. Some recent papers on the
subject like \cite{BaCiGaMa,BrWa,CiWa,Stolovitch-IHES,Stolovitch-Cartan} contain
some geometrical ingredients (commuting vector fields, symmetry groups), but the
methods used remain mostly analytical. In the present paper, we will follow the
geometrical approach in a more substantial way. We will use torus actions to
characterize Poincaré-Dulac normalization, and use this characterization to show
that any analytic vector field which is integrable in the non-Hamiltonian sense
admits a convergent Poincaré-Dulac normalization. More precisely, we have the
following

\begin{thm}
\label{thm:main} Let $X$ be a local analytic vector field in $(\bbK^n,0)$, where
$\bbK = \bbR$ or $\bbC$, with $X(0) = 0$. Suppose that there is a natural number
$m$, $1 \leq m \leq n$, such that there are $m$ local analytic vector fields $X_1 =
X, X_2, ..., X_m$ and $n-m$ local analytic functions $f_1,...,f_{n-m}$ in
$(\bbK^n,0)$ with the following properties :
\\
i) The vector fields $X_1 = X, X_2, ..., X_m$ commute pairwise, i.e.
\begin{equation}
[X_i,X_j] = 0 \ \ \forall \ i,j=1,...,m \ \ ,
\end{equation}
 and they are linearly independent almost everywhere,
i.e.
\begin{equation}
X_1 \wedge ... \wedge X_m \neq 0 \ .
\end{equation}
ii) The functions $f_1,...,f_{n-m}$ are common first integrals for $X_1,...,X_m$,
i.e.
\begin{equation}
X_i(f_j) = 0 \ \ \forall \ i=1,...,m \ ; \ j = 1,...,n-m \ \ ,
\end{equation}
 and they are functionally
independent almost everywhere, i.e.
\begin{equation}
df_1 \wedge ... \wedge df_{n-m} \neq 0 \ .
\end{equation}
Then there exists a local analytic Poincaré-Dulac normalization for $X$ in a
neighborhood of $0$ in $\bbK^n$. In other words, there is a local analytic system of
coordinates in $(\bbK^n,0)$ such that if we denote by $X^s$ the semi-simple part of
the linear part of $X$ with respect to this coordinate system then we have
\begin{equation}
[X, X^s ] = 0 \ .
\end{equation}
\end{thm}

{\it Remarks.}

1. In the above theorem, $m$ can be any natural number between 1 and $n$. When
$m=1$, we have a complete set of first integrals (and the vector field will be
automatically very resonant if $n>1$). When $m=n$ we have a complete set of
commuting infinitesimal symmetries. When $1 < m < n$, we have a mixture of commuting
infinitesimal symmetries and first integrals.

2. If a vector field satisfies the two conditions i) and ii) of the above theorem,
then one says that $X$ is integrable in the non-Hamiltonian sense (see e.g.
\cite{Bogoyavlenskij,ZungReduction}). The above theorem may viewed as the
non-Hamiltonian version of Theorem 1.1 of our previous paper \cite{ZungBirkhoff}.
Our proof of the above theorem is also based on a geometrical method developed in
\cite{ZungBirkhoff}.

3. When $n=2$ the above theorem has been obtained by Bruno and Walcher \cite{BrWa}.
Their proof is based on Bruno's results \cite{Bruno}, and therefore uses the fast
convergence method directly or indirectly, in contrast to our method.

4. It is interesting to compare the above theorem to similar but different results
of Stolovitch \cite{Stolovitch-IHES,Stolovitch-Cartan}. The similarity :
Stolovitch's results are also about the existence of a convergent normalization
under some integrability conditions. The difference : Stolovitch needs only formal
first integrals (it is roughly speaking the $A$-condition of Bruno, see
\cite{Bruno,Stolovitch-IHES}) while we use analytic (or eventually meromorphic - we
will consider the meromorphic case in a subsequent paper) ones, but to compensate
for this he also needs Diophantine conditions while we don't.

The reason why we are interested in integrability conditions while studying local
normal forms of vector fields is the following : both the search for first integrals
and symmetries and the search for a normalization are methods to reduce and solve
differential equations, and in good cases one can solve the equations completely by
both methods, so it is natural that the two methods should be very closely related.
Theorem \ref{thm:main}, as well as the title of the present paper, is about such
relations. (See \cite{ZungBirkhoff} for a more detailed discussion in the
Hamiltonian case).

The rest of this paper is organized as follows : In Section \ref{section:action} we
will use torus actions to characterize Poincaré-Dulac normalizations (see
Proposition \ref{prop:action}). As a corollary of this characterization, we obtain
that a real analytic vector field admits a local real analytic Poincaré-Dulac
normalization if and only if it admits a local holomorphic Poincaré-Dulac
normalization (see Proposition \ref{prop:realcomplex}). We will prove Theorem
\ref{thm:main} in Section \ref{section:proof}, using the above toric
characterization and a geometrical method of approximation developed in
\cite{ZungBirkhoff}. In Section \ref{section:isochore} we will extend our results to
the case of isochore (i.e. volume preserving) vector fields. In particular, we will
obtain an improvement of a result of Vey \cite{Vey2} about the existence of a
simultaneous normalization for an $(n-1)$-tuple of pairwise commuting holomorphic
isochore vector fields in $(\bbC^n,0)$ : we will show that, without any
nonresonance/nondegeneracy condition, Vey's theorem is still true, both in  $(\bbC^n,0)$
and in $(\bbR^n,0)$. \\

\section{Normalization and torus action}
\label{section:action}

In this section, we will recall some basic notions about Poincaré-Dulac normal
forms, and show that these normal forms are governed by torus actions.

Let $X$ be a given analytic vector field in a neighborhood of $0$ in $\bbK^n$, where
$\bbK = \bbR$ or $\bbC$, with $X(0) = 0$. When $\bbK = \bbR$, we may also view $X$
as a holomorphic (i.e. complex analytic) vector field by complexifying it. Denote by
\begin{equation}
X = X^{(1)} + X^{(2)} + X^{(3)} + ...
\end{equation}
the Taylor expansion of $X$ in some local system of coordinates, where $X^{(k)}$ is
a homogeneous vector field of degree $k$ for each $k \geq 1$. The algebra of linear
vector fields on ${\mathbb K}^{n}$, under the standard Lie bracket, is nothing but
the reductive algebra $gl(n, \bbK) = sl(n,\bbK) \oplus \bbK$. In particular, we have
\begin{equation}
X^{(1)} = X^s + X^{nil},
\end{equation}
where $X^s$ (resp., $X^{nil}$) denotes the semi-simple (resp., nilpotent) part of
$X^{(1)}$. There is a complex linear system of coordinates $(x_j)$ in ${\mathbb
C}^{n}$ which puts $X^s$ into diagonal form:
\begin{equation}
X^s = \sum_{j=1}^n \gamma_j x_j \partial / \partial x_j ,
\end{equation}
where $\gamma_j$ are complex coefficients, called {\it eigenvalues} of $X$ (or
$X^{(1)}$) at $0$.

For each natural number $k \geq 1$, the vector field $X^s$ acts linearly on the
space of homogeneous vector fields of degree $k$ by the Lie bracket, and the
monomial vector fields are the eigenvectors of this action:
\begin{equation}
[\sum_{j=1}^n \gamma_j x_j \partial / \partial x_j , x_1^{b_1}x_2^{b_2}...x_n^{b_n}
\partial / \partial x_l] = (\sum_{j=1}^n b_j\gamma_j - \gamma_l) x_1^{b_1}x_2^{b_2}...x_n^{b_n}
\partial / \partial x_l .
\end{equation}

When an equality of the type
\begin{equation}
\sum_{j=1}^n b_j\gamma_j - \gamma_l = 0
\end{equation}
holds for some nonnegative integer $n$-tuple $(b_j)$ with $\sum b_j \geq 2$, we will
say that the monomial vector field $x_1^{b_1}x_2^{b_2}...x_n^{b_n}
\partial / \partial x_l$ is a resonant term, and that the $n$-tuple $(b_1,...,b_l - 1,..., b_l)$
is a resonance relation for the eigenvalues $(\gamma_i)$. More precisely, a {\it
resonance relation} for  the $n$-tuple of eigenvalues $(\gamma_j)$ of a vector field
$X$ is an $n$-tuple $(c_j)$ of integers satisfying the relation $\sum c_j \gamma_j =
0,$ such that $c_j \geq -1, \sum c_j \geq 1,$ and at most one of the $c_j$ may be
negative. Notice here a major difference between the non-Hamiltonian case (which is
the case of this paper) and the Hamiltonian case (studied in \cite{ZungBirkhoff}) :
in the Hamiltonian case there are no positivity conditions like $c_j \geq -1$, so in
the Hamiltonian case the set of resonance relations for a given $n$-tuple of
frequencies is a sublattice in $\bbZ^n$, while in the non-Hamiltonian case the set
of resonance relations for a given $n$-tuple of eigenvalues is a convex set in
$\bbZ^n$ which does not contain any line.

Denote by
\begin{equation}
{\mathcal R} \subset {\mathbb Z}^n
\end{equation}
the subset of ${\mathbb Z}^n$ consisting of all resonance relations $(c_j)$ for a
given vector field $X$. The number
\begin{equation}
q = \dim_\bbZ (\rm linear \ hull \ of \ {\mathcal R} \ in \ {\mathbb Z}^n)
\end{equation}
is called the {\it degree of resonance} of $X$. Of course, the degree of resonance
depends only on the eigenvalues of the linear part of $X$, and does not depend on
the choice of local coordinates. If $q=0$ then we say that $X$ is {\it nonresonant}.

The vector field $X$ is said to be in {\it Poincaré-Dulac normal form} (or normal
form for short) if it commutes with the semisimple part of its linear part (see e.g.
\cite{Bruno,Roussarie}):
\begin{equation}
[X,X^s] = 0.
\end{equation}
The above equation means that if $X$ is in normal form then its nonlinear terms are
resonant. In particular, if $X$ is nonresonant then there are no resonant terms, and
$X$ is linear when it is in normal form. A transformation of coordinates which puts
$X$ in Poincaré-Dulac normal form is called a {\it Poincaré-Dulac normalization}. It
is a classical result of Poincaré and Dulac that any analytic vector field which
vanishes at 0 admits a {\it formal} Poincaré-Dulac normalization.

Denote by ${\mathcal Q} \subset {\mathbb Z}^n$ the integral sublattice of ${\mathbb
Z}^n$ consisting of $n$-dimensional vectors $(\rho_j) \in {\mathbb Z}^n$ which
satisfy the following properties :
\begin{equation}
\label{eqn:Q} \sum_{j=1}^n \rho_j c_j = 0 \  \forall \ (c_j) \in {\mathcal R} \ , \
{\rm and} \ \ \rho_j = \rho_k \ \ {\rm if} \ \ \gamma_j = \gamma_k \
\end{equation}
(where $\mathcal R$ is the set of resonance relations for the eigenvalues
$(\gamma_j)$ as before). We will call the number
\begin{equation}
\label{eqn:r} r = \dim_\bbZ \mathcal Q
\end{equation}
the {\it toric degree} of $X$ (or $X^{(1)}$). Of course, this number depends only on
the eigenvalues of the linear part of $X$, and we have the following inequality :
\begin{equation}
q + r \leq n \ ,
\end{equation}
where $q$ is the degree of resonance.

Let $(\rho^{1}_j),...,(\rho^r_j)$ be a basis of $\mathcal Q$. For each $k=1,...,r$
define the following diagonal linear vector field $Z_k$ :
\begin{equation}
\label{eqn:Z} Z_k = \sum_{j=1}^n \rho^k_j x_j \partial / \partial x_j \ .
\end{equation}
The vector fields $Z_1,...,Z_r$ have the following remarkable properties :

a) They commute pairwise and commute with $X^s$ and $X^{nil}$, and they are linearly
independent almost everywhere.

b) $iZ_j$ is a periodic vector field of period $2\pi$ for each $j \leq r$ (here $i =
\sqrt{-1}$). What does it mean is that if we write $iZ_j = \Re (iZ_j) + i \Im (iZ_j)
$, then $\Re (iZ_j)$ is a periodic real vector field in ${\mathbb C}^n = {\mathbb
R}^{2n}$ which preserves the complex structure.

c) Together, $iZ_1,..., iZ_r$ generate an effective linear ${\mathbb T}^r$-action in
${\mathbb C}^n$, which preserves $X^s$ and $X^{nil}$.

If $X$ is in Poincaré-Dulac normal form, i.e. $[X,X^s] = 0$, then we also have
\begin{equation}
[X,Z_k] = 0 \ \ \ \forall \  k=1,...,r.
\end{equation}

Indeed, if $W = x_1^{b_1}x_2^{b_2}...x_n^{b_n} \partial / \partial x_l$ is a
resonant monomial term, then we have $\sum b_j \gamma_j - \gamma_l = 0$, which
implies, via Equation (\ref{eqn:Q}), that $\sum_{j=1}^n b_j \rho^k_j - \rho_l = 0$,
and therefore $[Z_k, W] = (\sum_{j=1}^n b_j \rho^k_j - \rho_l) W = 0$. Hence if all
non-linear terms in $X$ are resonant then we have $[Z_k,X] = [Z_k,X^{(1)}] = 0$.

The above commutation relations mean that if $X$ is in normal form, then it is
preserved by the effective $r$-dimensional torus action generated by
$iZ_1,...,iZ_r$.

Conversely, suppose that there is a local analytic effective action of $\bbT^r$ in
$(\bbC^n,0)$ which preserves $X$, and whose linear part is generated by
$iZ_1,...,iZ_r$ (in other words, its weights are given by the lattice $\mathcal Q$).
Then using the classical Bochner's linearization theorem for compact group actions
(via the standard averaging method), we can linearize this torus action, i.e. we may
suppose that the action is linear (in an appropriate local analytic coordinate
system). Then this action is actually generated by $iZ_1,...,iZ_r$. Since the action
preserves $X$ by assumptions, we have that $[X,Z_1] = ... = [X,Z_r] = 0$. But by
definition of $Z_1,...,Z_r$, the semisimple part of the linear part of $X$ is a
linear combination of these vector fields $Z_1,...,Z_r$, hence we have $[X,X^s] =
0$. In other words, a linearization of our torus action is also a Poincaré-Dulac
normalization for $X$. Thus we have proved the following result :

\begin{prop}
\label{prop:action} A holomorphic vector field $X$ in a neighborhood of $0$ in
${\mathbb C}^n$ admits a convergent Poincaré-Dulac normalization if and only if it
is preserved by an effective holomorphic action of a real torus of dimension $r$,
where $r$ is the toric degree of $X^{(1)}$ as defined in (\ref{eqn:r}), in a
neighborhood of $0$ in ${\mathbb C}^n$, which has $0$ as a fixed point and whose
linear part at $0$ has appropriate weights (given by the lattice $\mathcal Q$
defined in (\ref{eqn:Q}), which depends only on the linear part $X^{(1)}$ of $X$).
$\diamondsuit$
\end{prop}

{\it Remark}. The above proposition is true in the formal category as well. But of
course, any vector field admits a formal Poincaré-Dulac normalization, and a formal
torus action.

When $X$ is real analytic, the torus action in the above proposition still lives in
the complex space in general. But in this case, there is a natural complex
conjugation in the torus action, which allows us to deduce from Proposition
\ref{prop:action} the following result :

\begin{prop}
\label{prop:realcomplex} A real analytic vector field $X$ in a neighborhood of $0$
in ${\mathbb R}^n$ with $X(0) = 0$ admits a local real analytic Poincaré-Dulac
normalization if and only if it admits a local holomorphic Poincaré-Dulac
normalization when considered as a holomorphic vector field.
\end{prop}

The proof of Proposition \ref{prop:realcomplex} is absolutely similar to the proof
of Proposition 1.3 of our previous paper \cite{ZungBirkhoff} (its Hamiltonian
version), so we will omit it here. $\diamondsuit$

\section{Torus actions for integrable vector fields}
\label{section:proof}

{\it Proof of Theorem  \ref{thm:main}}. Invoking Proposition \ref{prop:action}, we
will prove Theorem \ref{thm:main} by finding a torus action. In view of Proposition
\ref{prop:realcomplex}, it suffices to prove Theorem \ref{thm:main} in the complex
case, so in this section we will assume that $\bbK = \bbC$.

Recall that a vector field $X$ is said to be in normal form up to order $M$ (where
$M$ is a natural number) in a given coordinate system if $[X,X^s] = O(|x|^M)$ (i.e.
$[X,X^s]$ does not contain terms of order $< M$), where $X^s$ denotes the semisimple
part of the linear part of $X$ in this coordinate system. A coordinate
transformation that puts $X$ in normal form up to order $M$ is called a
normalization up to order $M$. Such a local holomorphic normalization up to order
$M$ always exists (for any $M$), according to the classical theorem of Poincaré and
Dulac. We have the following lemma, whose proof is straightforward and uses only
elementary linear algebra :

\begin{lem}
\label{lemma:orderM} Suppose that a vector field $X$ is in Poincaré-Dulac (formal or
holomorphic) normal form. Denote by $r$ the toric degree of $X$, and by
$Z_1,...,Z_r$ the vector fields defined by Equation (\ref{eqn:Z}). If $Y$ is a
vector field which commutes with $X$ then $[Z_k,Y] = 0$ for $k = 1,...,r$. If $f$ is
a first integral for $X$, i.e. $X(f) = 0$, then $Z_k(f) = 0$ for $k = 1,...,r$. If
$X$ is in normal form up to order $M$ (for some natural number $M$) then $[Z_k,Y] =
O(|x|^M)$ and $Z_k(f) = O(|x|^M)$ for $k=1,...,r$. $\diamondsuit$
\end{lem}

Now let $X$ be a local holomorphic vector field in $(\bbC^n,0)$ which satisfies the
conditions of Theorem \ref{thm:main}, with $X_1 = X,...,X_m$ and $f_1,...,f_{n-m}$
being the $m$ local holomorphic pairwise commuting vector fields and $(n-m)$
holomorphic common first integrals respectively. By assumptions, $X_1 \wedge ...
\wedge X_m \neq 0$ and $df_1 \wedge ... \wedge df_{n-m} \neq 0$ almost everywhere.

Fix a holomorphic coordinate system $x = (x_j)$ in ${\mathbb C}^n$, a standard
Hermitian metric in ${\mathbb C}^n$ which goes with it, and a sufficiently small
positive number $\epsilon_0$. Denote by
\begin{equation}
S = \{ x \in {\mathbb C}^n \ , \ |x| < \epsilon_0 \ , \ X_1 \wedge X_2 \wedge ...
\wedge X_m = 0 \ or \ df_1 \wedge ... \wedge df_{n-m}(x) = 0 \}
\end{equation}
the singular locus of our $n$-tuple of vector fields and functions
$X_1,...,X_m,f_1,...,f_{n-m}$. By assumptions, $S$ is a complex analytic set of
complex codimension at least 1. In particular, we have the following \L
ojasiewicz-type inequalities (see e.g. \cite{Lojasiewicz}): there exist a natural
number $N$ and a positive constant $C$ such that
\begin{equation}
|X_1 \wedge ... \wedge X_m (x) | \geq C (d(x, S))^N
\end{equation}
and
\begin{equation}
|df_1 \wedge ... \wedge df_{n-m} (x) | \geq C (d(x, S))^N
\end{equation}
for any $x$ with $|x| < \epsilon_0$, where the norms applied to $X_1 \wedge ...
\wedge X_m (x)$ and $df_1 \wedge ... \wedge df_{n-m} (x)$ are some standard norms in
the space of $m$-vectors and $(n-m)$-vectors respectively, and $d(x,S)$ is the
distance from $x$ to $S$ with respect to the Euclidean metric.

For each $d \in {\mathbb N}$ and a small positive number $\epsilon (d) > 0$ (which
will be chosen later in function of $d$, with $\lim_{d \to \infty} \epsilon (d) =
0$), define the following open subset $U_{d,\epsilon (d)}$ of ${\mathbb C}^n$:
\begin{equation}
U_{d,\epsilon (d)} = \{ x \in {\mathbb C}^n \ , \ |x| < \epsilon (d) \ , \ d(x,S) >
|x|^d \} \ .
\end{equation}
For each $k=1,...,r$, we will define a holomorphic vector field $\cZ_k$ in
$U_{d,\epsilon(d)}$, such that (the real part of) $\sqrt{-1}\cZ_k$ is a  periodic
vector field of period $2\pi$, and in such a way that for any two natural numbers
$d_1 \neq d_2$ the vector field $\cZ_k$ defined for $U_{d_1,\epsilon(d_1)}$
coincides with the vector field (with the same name) $\cZ_k$ defined for
$U_{d_2,\epsilon(d_2)}$ on the intersection $U_{d_1,\epsilon(d_1)} \cap
U_{d_1,\epsilon(d_1)}$, as follows.

Let
\begin{equation}
(x^d_j) = \Phi^{D(d)}(x_j) = (x_j) + \ {\rm higher \ order \ terms}
\end{equation}
be a local holomorphic coordinate transformation which puts $X$ into normal form up
to order $D(d)$, where $D(d)$ is a sufficiently large number (as large as we wish),
which will be chosen in function of $d$, with $\lim_{d \to \infty} D(d) = \infty$.
We can choose $\Phi^{D(d)}$ in such a way that for any two natural numbers $d_1 \neq
d_2$ we have
\begin{equation}
\Phi^{D(d_1)}(x_j) = \Phi^{D(d_2)}(x_j) + \ {\rm terms \ of \ order} \geq \min
(D(d_1), D(d_2))
\end{equation}
Denote by
\begin{equation}
Z^d_k = \sum_{j=1}^n \rho^k_j x^d_j \partial / \partial x^d_j
\end{equation}
the vector fields that are defined as in Equation (\ref{eqn:Z}), but for the
semisimple part of $X$ with respect to the local coordinate system $(x^d_j)$. In
particular, $\sqrt{-1}Z_k$ is a holomorphic periodic vector field of period $2\pi$
for each $k=1,...,r$. According to Lemma \ref{lemma:orderM}, we have
\begin{equation}
[Z^d_k,X_j] (x) = O(|x|^{D(d)}) \ \ \forall j=1,...,m,
\end{equation}
and
\begin{equation}
Z^d_j (\boldf) (x) = O(|x|^{D(d)})
\end{equation}
where $\boldf = (f_1,...,f_{n-m})$ denotes the vector-valued common first integral
of $X_1,...,X_m$. In other words, $\boldf$ is a vector-valued first integral of
$Z^d_k$ up to order $D(d)$, and $X_1,...,X_m$ commute with $Z^d_j$ up to order
$D(d)$.

Let $y$ be an arbitrary point in $U_{d,\epsilon (d)}$. Then, due to the \L
ojasiewicz-type inequalities and the definition of $U_{d,\epsilon (d)}$, we have:
\begin{equation}
|X_1 \wedge ... \wedge X_m (y)| > C |y|^{dN}
\end{equation}
and
\begin{equation}
|df_1 \wedge ... \wedge df_{n-m}(y)| > C |y|^{dN} .
\end{equation}

Denote by $\Gamma^d_k(t) = \Gamma^d_k(t,y)$ the closed curve ($t \in [0,2\pi]$)
which is the orbit of the periodic  vector field $\Re (\sqrt{-1}Z^d_k)$ which begins
at $y$. We have that $\Gamma^d_k(0)= y$, and $\frac{1}{2}|y| \leq |\Gamma^d_k(t)|
\leq 2|y|$ for any $t \in [0,2\pi]$, provided that $\epsilon (d)$ is small enough.

It follows from the fact that $\boldf$ is a vector-valued first integral of $Z^d_k$
up to order $D(d)$ and $X_1,...,X_m$ commute with $Z^d_k$ up to order $D(d)$, that
we have the following inequalities, provided that $D(d)$ is large enough (say $D(d)
> 4dN$) and $\epsilon (d)$ small enough:
\begin{equation}
\label{eqn:4}
\begin{array}{l}
|\boldf(z) - \boldf(y)| < |y|^{D_1(d)} \\
 |[X_j,Z^d_k](z)| < |y|^{D_1(d)} \ \ \forall j=1,...,m \\
|X_1 \wedge ... \wedge X_m (z)| > \frac{C}{2} |y|^{dN}  \\
|df_1 \wedge ... \wedge df_{n-m}(z)| > \frac{C}{2} |y|^{dN}
\end{array}
\end{equation}
for any point $z$ lying on the curve $\Gamma^d_k$. Here $D_1(d)$ is some
sufficiently large number to be chosen in function of $d$. (We may choose $D_1(d)$
as large as we need, and then choose $D(d)$ and $\epsilon (d)$ correspondingly so
that the above inequalities hold).

The above inequalities imply the following things :

a) The curve $\Gamma^{d}_k$ is very close to the (regular part of the) level set
$L_y = \boldf^{-1}(\boldf(y))$ of the vector-valued first integral
$\boldf=(f_1,...,f_{n-m})$ which contains the point $y$, in the sense that it can be
projected orthogonally to a smooth closed curve $\hat{\Gamma}^d_k(t)$ ($t \in [1,
2\pi]$) lying on $L_y$ which is very close to $\Gamma^d_k$ in $C^1$-topology : the
distance from $\hat{\Gamma}^d_k$ to $\Gamma^d_k$ in $C^1$-topology is bounded from
above by $|y|^{D_2(d)}$, where $D_2(d)$ is a sufficiently large number depending on
$d$. (We may take $D_2(d)$ as large as we want, and then choose $D_1(d)$ and
$\epsilon (d)$ correspondingly so that the above upper bound holds true).

b) The regular part of the level set $L_y$ is of complex dimension $m$, and its
tangent space at each point is spanned by the vectors $X_1,...,X_m$. Moreover, the
regular part of $L_y$ has a flat affine structure given by the vector fields
$X_1,...,X_m$, because these vector fields commute.

c) If we write  $d\hat{\Gamma}^d_k(t)/dt = \sum_{j=1}^m \Re (a^j_k (t) X_j
(\hat{\Gamma}^d_k(t)))$, then the complex functions $a^j_k (t)$ are nearly constant
in the sense that $|a^j_k (t) - a^j_k (0)| \leq |y|^{D_3(d)}$ for $t \in [1, 2\pi]$,
where $D_3(d)$ is a sufficiently large number depending on $d$ (we may take $D_3(d)
= D_2(d) -1$). This fact follows from the almost commutativity of the vector fields
$X_1,...,X_m$ with the vector field $Z^d_k$ (the second line of (\ref{eqn:4})), and
the fact that $|d\hat{\Gamma}^d_k(t)/dt - \Re (\sqrt{-1}Z^d_k(\hat{\Gamma}^d_k(t)))|
< |y|^{D_2(d)}$ by the above point a).

d) By approximation (implicit function theorem), there exist complex numbers
$a^1_k,...,a^m_k$ such that $|a^j_k - a^j_k (0)| < |y|^{D_3(d)}$, and that the
time-$2\pi$ flow of the vector field $\sum_{j=1}^m a^j_k X_j$ on $L_y$ fixes point
$y$. In other words, the real vector field $\Re (\sum a^j_k X_j)$ has a periodic
orbit of period $2\pi$ which passes via $y$, and this orbit is $C^1$-close to
$\hat{\Gamma}^d_k(t,y)$.

e) Due to the flat affine structure of $L_y$, the numbers $a^1_k,...,a^m_k$ are well
defined, i.e. unique. And they don't depend on the choice of $y$ in $L_y$ (at least
locally). We may consider $a^1_k,...,a^m_k$ as functions of $y$ :
$a^1_k(y),...,a^m_k(y)$. These are holomorphic functions (due to the holomorphic
implicit-function theorem) which are constant on the connected components in
$U_{d,\epsilon (d)}$ of the level sets of the vector-valued function $\boldf$, and
are uniformly bounded in $U_{d,\epsilon (d)}$ by 1 (provided that $\epsilon (d)$ is
small enough).

Now define the vector field $\cZ_k$ as follows :
\begin{equation}
\cZ_k (y) = -\sqrt{-1} \sum_{j=1}^m a^j_k (y) X_j (y)
\end{equation}

Then $\cZ_k$ is a holomorphic vector field in $U_{d,\epsilon (d)}$ with the
following remarkable properties (for each $k=1,...,r$):

a) $\cZ_k$ is uniformly bounded by 1, and $\sqrt{-1}\cZ_k$ is a periodic vector
field of period $2\pi$ (at least in some open subset of $U_{d,\epsilon (d)}$).

b) If $\cZ_k$ is a vector field  defined as above for $U_{d, \epsilon (d)}$, and
$\cZ'_k$ is also a vector field defined as above but for $U_{d', \epsilon (d')}$
with $d' \neq d$, then $\cZ_k$ coincides with $\cZ'_k$ on the intersection of $U_{d,
\epsilon (d)}$ with $U_{d', \epsilon (d')}$, due to the uniqueness of $a^j_k$.
Indeed, the vector field $\cZ_k$ commutes with the vector field $\cZ'_k$ on $U_{d,
\epsilon (d)} \cap U_{d', \epsilon (d')}$ by construction, their difference $\cZ_k -
\cZ'_k$ is tangent to the level sets of $\boldf$ in $U_{d, \epsilon (d)} \cap U_{d',
\epsilon (d')}$ and is a constant vector field with respect to the flat affine
structure on each level set, and $\sqrt{-1}(\cZ_k - \cZ'_k)$ is periodic of period
at most $2\pi$ there. But the coefficients of $\cZ_k (y) - \cZ'_k (y)$ when written
as a linear combination of $X_1,...,X_m$ are bounded from above by $|y|^{min
(D_3(d),D_3(d'))}$, i.e. $\sqrt{-1}(\cZ_k - \cZ'_k)$ is too small to be periodic of
period $2\pi$ unless it is zero. Thus we have $\cZ_k = \cZ'_k$ in $U_{d, \epsilon
(d)} \cap U_{d', \epsilon (d')}$.

In other words, we have defined $r$ bounded holomorphic vector fields
$\cZ_1,...,\cZ_r$ on $U = \bigcup_{d=1}^{\infty} U_{d, \epsilon (d)}$. These vector
fields commute pairwise due to the flat structure of the level sets $L_y$ (they are
constant on each $L_y$ with respect to the flat structure). And
$\sqrt{-1}\cZ_1,...,\sqrt{-1}\cZ_r$ are periodic of period $2\pi$ (ar least in some
open subset). The following lemma taken from \cite{ZungBirkhoff} shows that
$\cZ_1,...,\cZ_r$ can be extended holomorphically in a neighborhood of $0$ in
${\mathbb C}^n$, i.e. there are holomorphic vector fields in a neighborhood of $0$
in ${\mathbb C}^n$ which coincide with them on $U$. Thus we found the generators for
an effective torus action of dimension $r$ which preserves $X= X_1$ (and
$X_2,...,X_m$ as well). Due to our approximation process, the linear part of $\cZ_k$
is the same as the linear part of $Z^d_k$ (for any $d$), i.e. the weights of our
torus action are appropriate in the sense of proposition \ref{prop:action}. Thus we
can apply Proposition \ref{prop:action} to finish the proof of Theorem
\ref{thm:main}. $\diamondsuit$

\begin{lem}[\cite{ZungBirkhoff}]
Let $U = \bigcup_{d=1}^{\infty} U_d$, with $U_d = \{ x \in {\mathbb C}^n, |x| <
\epsilon_0, d(x,S) > |x|^d \}$, where $\epsilon_d$ is an arbitrary series of
positive numbers and $S$ is a local proper complex analytic subset of ${\mathbb
C}^n$  ($codim_{\mathbb C} S \geq 1$). Then any bounded holomorphic function in $U$
has a holomorphic extension in a neighborhood of $0$ in ${\mathbb C}^n$.
\end{lem}

The proof of the above lemma is given in \cite{ZungBirkhoff}, so we will not repeat
it here. $\diamondsuit$

As a corollary of Theorem \ref{thm:main}, we get the following result about
simultaneous Poincaré-Dulac normal forms for commuting vector fields :

\begin{thm}
\label{thm:main2} Any $m$-tuple of pairwise commuting analytic vector fields
$X_1,...,X_m$ in a neighborhood of $0$ in $\bbK^n$ ($n \geq m \geq 1$), which are
linear independent almost everywhere, and which have $(n-m)$ functionally
independent common analytic first integrals, admits a simultaneous convergent
Poincaré-Dulac normalization.
\end{thm}

{\it Proof}. By assumptions, each of the vector fields $X_1,...,X_m$ is integrable
in the non-Hamiltonian sense, so the proof of Theorem \ref{thm:main} shows that for
each of them there is a torus action whose linearization normalizes it. Moreover, by
construction, all these torus actions commute with each other and preserve all the
vector fields $X_1,...,X_m$. Thus we can combine these torus actions to get one
``big'' torus action whose linearization will normalize $X_1,...,X_m$
simultaneously. $\diamondsuit$ \\

\section{The isochore case}
\label{section:isochore}

A normalization of an isochore (i.e. volume-preserving) vector field in $({\mathbb
C}^n, 0)$ is a Poincaré-Dulac normalization $\tilde{x} = \phi(x)$ which preserves
the volume form $dx_1 \wedge ... \wedge dx_n$. It is not surprising that an isochore
vector field always admits a formal normalization.  Notice that in the isochore
case, there is at least one resonance relation : the sum of the eigenvalues is zero.
This relation implies that the  vector fields $Z_1,...,Z_r$ defined in (\ref{eqn:Z})
are also isochore. Naturally, Proposition \ref{prop:action} also has an isochore
version.

\begin{prop}
\label{prop:isochoreaction} A holomorphic vector field $X$ in a neighborhood of $0$
in ${\mathbb C}^n$ admits a convergent Poincaré-Dulac normalization if and only if
it is preserved by an effective holomorphic isochore action of a real torus
${\mathbb T}^{r^I}$ a neighborhood of $0$ in ${\mathbb C}^n$, where $r$ is the toric
degree of $X$, which has $0$ as a fixed point and whose linear part at $0$ has
appropriate weights (given by the lattice ${\mathcal Q}$ defined in (\ref{eqn:Q}),
which depends only on the linear part $X^{(1)}$ of $X$).
\end{prop}

The proof of Proposition \ref{prop:isochoreaction} is absolutely similar to the
proof of Proposition \ref{prop:action}. $\diamondsuit$

Similarly, Proposition \ref{prop:realcomplex} remains true in the isochore case as
well.

In \cite{Vey2}, Vey showed that if $n-1$ pairwise commuting isochore vector fields
in ${\mathbb C}^n$ have linearly independent diagonalizable linear parts, then they
are simultaneously normalizable.  Though Vey's proof didn't use the fast convergence
method directly, it used a result of Malgrange on ``Frobenius with singularities''
\cite{Malgrange}, which in turn had been proved by the fast convergence method. We
will show that Vey's theorem is still true without any nondegeneracy condition :

\begin{thm}
\label{thm:isochore} Any $(n-1)$-tuple of pairwise commuting analytic isochore
vector fields in a neighborhood of $0$ in ${\mathbb K}^n$, which are linearly
independent almost everywhere, admits a simultaneous convergent normalization.
\end{thm}

Note that in the above theorem, though the vector fields are assumed to be linearly
independent almost everywhere, their linear parts may be very degenerate and may
indeed be linearly dependent.

We will only give here a sketch of the proof of the above theorem. It consists of
the following steps:

1) Find a first integral: The 1-form $\alpha = i_{X_1}i_{X_2}...i_{X_{n-1}} \Omega$,
where $\Omega$ is the volume form which is preserved by $(n-1)$ pairwise commuting
vector fields $X_1,...,X_{n-1}$, is closed, hence exact (Poincaré's lemma). If $g$
is a function such that $d g = \alpha$ then $g$ is a common first integral for our
vector fields.

2) Use Theorem \ref{thm:main2}, since now we have $n-1$ commuting vector fields and
1 first integral.

3) Keep track of the isochore condition. (Everything is made in an isochore way).

Details are left to the reader. \\

{\bf Acknowledgements}. I would like to thank the referee for his critical remarks.
\\

\bibliographystyle{amsplain}

\end{document}